\documentclass[11pt]{article}

\usepackage[utf8]{inputenc}
\usepackage[english]{babel}
\usepackage[all]{xy}
\usepackage{tikz}
\usepackage{hyperref}
\usepackage{amsmath}
\usepackage{amssymb}
\usepackage{amsthm}
\usepackage{mathrsfs}
\usepackage{mathtools}
\usepackage{amsfonts}

\usepackage[a4paper, top=2.5cm, bottom=2.5cm, left=2.5cm, right=2.5cm]{geometry}

\usepackage[draft, textsize=tiny]{todonotes}

\newcommand{\del}{\partial}

\newcommand{\N}{\mathbb{N}}
\newcommand{\C}{\mathbb{C}}
\newcommand{\R}{\mathbb{R}}

\newtheorem{dfn}{Definition}
\newtheorem{exa}{Example}

\newtheorem{lem}{Lemma}
\newtheorem{prp}{Proposition}
\newtheorem{thm}{Theorem}
\newtheorem{cor}{Corollary}

\usetikzlibrary{matrix}

\begin{document}

    \title{A class of globally analytic hypoelliptic operators on compact Lie groups}

    \author{
    Max Reinhold Jahnke
    \footnote{
        The first author was supported by the Deutsche Forschungsgemeinschaft (DFG), grant JA 3453/1-1.}
    \\
    \href{mailto:max.jahnke@uni-koeln.de}{max.jahnke@uni-koeln.de} \\
    \\
    Mathematics Institute, \\
    University of Cologne
    \and Nicholas Braun Rodrigues \footnote{
        The second author was partially supported by Fundação de Amparo à pesquisa do Estado de São Paulo (FAPESP), grant 2020/15368-7.} 
    \\\href{mailto:braun@ime.unicamp.br}{braun@ime.unicamp.br} \\
    \\
    Department of Mathematics, \\
    State University of Campinas}
    
    \maketitle
    
    \begin{abstract}
        We obtain global analytic hypoellipticity for a class of differential operators that can be expressed as a zero-order perturbation of a sum of squares of vector fields with real-analytic coefficients on compact Lie groups. The key conditions are: the vector fields must satisfy Hörmander's finite type condition; there exists a closed subgroup whose action leaves the vector fields invariant; and the operator must be elliptic in directions transversal to the action of the subgroup. This paves the way for further studies on the regularity of sums of squares on principal fiber bundles.
    \end{abstract}
    
    \section{Introduction}
    \label{sec:introduction}

    We proved that zero-order perturbations of a certain class of differential operators expressed as a sum of squares of vector fields with real coefficients on compact Lie groups is globally analytic hypoelliptic. Our work generalizes the global analytic hypoellipticity results in \cite{braun_rodrigues_lower_2016} to the setting of compact Lie groups. Unlike \cite{braun_rodrigues_lower_2016}, we do not assume that the ambient manifold is a product of manifolds in order to separate ellipticity in one variable and work to obtain the regularity in the other ones. We instead use the intrinsic geometrical properties of the Lie groups to obtain the global regularity result. We assume three conditions. Naturally, we assume that the vector fields must satisfy Hörmander's finite type condition. To leverage the geometrical properties of the Lie group, we assume that there exists a closed subgroup whose action leaves the vector fields invariant and that we have ellipticity in directions transversal to the action of the subgroup. We believe that our results open doors to further investigation of the regularity of certain operators on principal fiber bundles. 
    
    The advantage of working first on Lie groups is that it is possible to give a global expression for the operators. The global expression makes it easy to discuss the commutation problem of the operator under study with the Laplacian, which would be more complicated in the case of fiber bundles in general. Furthermore, even though it is common to use even-order Sobolev spaces when using elliptic theory, we need to consider the fractional case, which we obtain by using interpolation theory. In order to keep the text self-contained, we include the results regarding interpolation in the appendix.

    Although the initial motivation for our research comes from the original work of Cordaro and Himonas \cite{cordaro_global_1994, cordaro_global_1998}, we were motivated by several other papers in this and related topics, such as \cite{silva_globally_2023, chinni_global_2017, chinni_microlocal_2022, barostichi_global_2017, bove_analytic_2018, ferra_global_2020, ferra_perturbations_2023, de_lessa_victor_fourier_2021, furdos_theorem_2022,araujo_global_2022, alexandrov_vanishing_2001, himonas_analyticity_2020}. The challenge of understanding the conditions that guarantee global analytic hypoellipticity of sum-of-squares operators on compact real-analytic manifolds remains an open and actively pursued area of research.

    Let $G$ be a compact Lie group with Lie algebra $\mathfrak g$ and let $K \subset G$ be a closed subgroup with Lie algebra $\mathfrak k$. We consider $X_0, X_1, \dots, X_\nu$, real and real-analytic $K$-invariant vector fields on $G$ satisfying the \emph{Hörmander finite-type condition} (H), that is, the vector fields, together with their brackets of length at most $r$, span the tangent space at every point of $G$, and such that
        \begin{equation}
        \label{eq:e1}
        \tag{E1}
            \operatorname{span} \{ \pi_\ast(X_1(g)), \dots, \pi_\ast(X_{\nu}(g)) \} = T_{\pi(g)} G/K, \forall g \in G,
        \end{equation} with $\pi : G \longrightarrow G/K$ the projection. 

        We recall that an operator $P$ is called \emph{globally analytic hypoelliptic} (GAH) on $G$ if for all distributions $u \in \mathcal{D}'(G)$ such that $Pu \in C^\infty(G)$ we have that $u \in C^\infty(G)$. 
        
        The following is our main result:
        
        \begin{thm}
        \label{thm:global}
            The operator
            \begin{equation}
            \label{eq:P}
                P = \sum_{j=1}^\nu X_j^2 + X_0 + a,
            \end{equation} with 
            $X_0, X_1, \dots, X_\nu$ real and real-analytic $K$-invariant vector fields on $G$ 
            satisfying conditions (H) and (E1), is (GAH) on $G$ for all $a \in C^\omega(G)$.
        \end{thm}

        Notice that, because of the zero-order term $a$, which does not need to be invariant under the action of $K$, the operator \eqref{eq:P} does not commute with the Laplace-Beltrami operator on $K$ coming from an ad-invariant metric.  The proof of the theorem uses the geometrical structure of the group $G$, together with the following local-global version of the theorem. Let $K$ be a compact Lie group with Lie algebra $\mathfrak k$, $U \subset \R^m$ be an open set and $Z_1, \ldots, Z_n \in \mathfrak{k}$. We denote the coordinates on $U$ by $(t_1, \ldots, t_m)$ and we consider the vector fields
        \begin{equation}
        \label{eq:vf2}
            L_j = \sum_{l=1}^n a_{jl}(t) \frac{\del}{\del t_l} + \sum_{l=1}^m b_{jl}(t) Z_l, \quad j = 1, \ldots, \nu, 
        \end{equation} defined on $U \times K$, with $a_{jl}, b_{jl}$ real-analytic functions on $U$ and satisfying the Hörmander condition and ellipticity in $t$: 

        \begin{equation} \label{eq:e2} \tag{E2}
            \operatorname{span} \left\{ \sum_{l=1}^n a_{jl}(t) \frac{\del}{\del t_l}: ~j = 1, \ldots, \nu \right\} = T_t(U), \qquad \forall t \in U.
        \end{equation}

        \begin{thm}
        \label{thm:local_global_main}
            The operator
            \begin{equation}
            \label{eq:Q}
                Q = \sum_{j=1}^\nu L_j^2 + L_0 + a,
            \end{equation} with $L_j$ as in \eqref{eq:vf2} satisfying (H) and (E2) is (GAH) on $U \times K$ for all $a \in C^\omega(U \times K)$.
        \end{thm}

    \begin{exa}
        Recall that the $\operatorname{SU}(2)$ is defined by
        $$ \operatorname{SU}(2) \doteq \left\{ \left(  \begin{matrix} z_1 & -\overline{z}_2 \\ z_2 & \overline{z}_1 \end{matrix} \right) : z_1, z_2 \in \C, ~|z_1|^2 + |z_2|^2 = 1 \right\}$$ with Lie algebra denoted by $\mathfrak{su}(2)$ and generated by
        $$ X = \left(  \begin{matrix} 0 & i \\ i & 0 \end{matrix} \right), \quad Y = \left(  \begin{matrix} 0 & -1 \\ 1 & 0 \end{matrix} \right), \quad  T = \left(  \begin{matrix} i & 0 \\ 0 & -i \end{matrix} \right).$$

        The vectors are invariant under the right-action of $G$ on itself, thus it is also invariant under the action of $K \doteq \exp(\operatorname{span}_{\R}\{T\})$. The following relation between the brackets of $X, Y$ and $T$
        \begin{equation}
            \label{sec:su2:basic_commutator_relations}
            [T, X] = 2Y, \quad [T, Y] = -2X, \quad [X, Y] = 2T
        \end{equation}
        guarantee that the operator $P = (fX)^2 + (gY)^2 + a$, with $a, f, g \in C^\omega(G)$ such that $f,g$ are non-vanishing functions that are invariant under the action of $K$, satisfy condition (H) and clearly also satisfy $(E1)$.
    \end{exa}

    The proof of the theorems uses two fundamental results, one in the theory of regularity of partial differential equations and another about the structure of Lie groups. In the following, we state both results, starting by the regularity result. First, we need a few definitions.

    Let $L = \{L_k\}_{k \in \N}$ be an increasing sequence of positive numbers such that $L_0 = 1$, $k \leq L_k$ and $L_{k+1} \leq CL_k$ for a constant $C > 0$. If $\Omega \subset \R^n$ is an open set, we denote by $C^L(\Omega)$ the set of functions $u$ in $\Omega$ such that, for every compact $K \subset \Omega$, there exists a constant $C_K > 0$ such that $$ \| D^\alpha u \|_{L^2(K)} \leq C_K (C_K L_{|\alpha|})^{|\alpha|}$$ for any multi-index $\alpha$. Notice that when $L_k = k+1$, the the space $C^L(\Omega)$ is the space of real-analytic functions in $\Omega$. In this case, we denote $L$ by $\omega$. 

    \begin{dfn}
        Let $x_0 \in \Omega$, $\xi_0 \in \R^n \backslash \{0\}$, $u \in \mathcal D'(\Omega)$ and $P$ be a differential operator of order $m$ with coefficients in $\Omega$; we say that $(x_0, \xi_0)$ is in the complement of the wavefront $WF_L(u; P)$ of $u$ with respect to $C^L$ and the iterates of $P$ if and only if there exists an open neighborhood $U$ of $x_0$, an open conic neighborhood $\Gamma$ of $\xi_0$ and a sequence $f_N \in \mathcal E'(\Omega)$ such that, for $f_N \doteq P^Nu$ in $U$, we have:
        $$ |\hat f_N (\xi)| \leq C(C(L_{mN} + |\xi|))^{nN + M}\qquad \xi \in \R^n, $$
        $$ |\hat f_N (\xi)| \leq C(CL_{mN})^{mN}(1 + |\xi|)^m, \quad \xi \in \Gamma, N = 0, 1, \ldots$$
        for constants $C > 0$ and $M$.
    \end{dfn}

    \begin{thm}[Bolley-Camus-Mattera (1979)]
        Let $u$ in $\mathcal D'(\Omega)$ and $P$ be a differential operator of order m with analytical coefficients in $\Omega$. Then: $$ \operatorname{WF}_L(u) \subset \operatorname{WF}_L(u; P) \cup \{ (x, \xi) \in T^*(\Omega)\backslash \{0\}; P_m(x, \xi) = 0\}$$
        in which $P_m$ being the principal part of $P$.
    \end{thm}
    
    \section{Preliminaries}\label{sec:preliminaries}

    In this section, we briefly recall facts about to work with right invariant vector fields and how to use Fourier analysis on compact Lie groups. Let $G$ be a compact Lie group with Lie algebra $\mathfrak g$. We denote $\mathfrak{X}_L(G)$ and $\mathfrak{X}_R(G)$ the set of all left and right invariant vector fields on $G$, \textit{i.e.} denoting by $L_h$ and $R_h$ the left and right translations by $h \in G$, a vector field $\chi$ belongs to $\mathfrak{X}_L(G)$ if $(L_h)_\ast \chi|_x = \chi|_{hx}$ for every $x,h \in G$, and $\chi$ belongs to $\mathfrak{X}_R(G)$ if $(R_h)_\ast \chi|_x = \chi|_{xh}$ for every $x,h \in G$. If $X \in \mathfrak{g}$, we denote by $X^L$ and $X^R$ the left and right invariant vector fields given by
    $$ X^L f(x) \doteq \frac{\mathrm{d}}{\mathrm{d} t}\Big|_{t=0} f(x \exp(tX)) \quad \text{and} \quad X^R f(x) \doteq \frac{\mathrm{d}}{\mathrm{d} t}\Big|_{t=0} f(\exp(tX) x), $$
    where $\exp$ stands for the exponential function from $\mathfrak{g}$ to $G$. The application that maps $\mathfrak{g}$ to $\mathfrak{X}_L$ and $\mathfrak{X}_G$ is a Lie algebra isomorphism.
    Let $\langle \cdot, \cdot \rangle_\mathfrak{g}$ be an Ad-invariant inner product on $\mathfrak{g}$, \textit{i.e.} $$\langle \text{Ad}(x) X , \text{Ad}(x) Y \rangle_\mathfrak{g} =  \langle  X ,  Y \rangle_\mathfrak{g}, \quad \forall x \in G, \forall X,Y \in \mathfrak{g},$$ where $\text{Ad}(x) = (L_x\circ R_{x^{-1}})_\ast$. In particular, $\langle \cdot, \cdot \rangle_\mathfrak{g}$ is also ad-invariant, in other words,
    $$\langle  [X,Y] ,  Z \rangle_\mathfrak{g} = \langle  X ,  [Y,Z] \rangle_\mathfrak{g} ,\quad \forall X,Y,Z \in \mathfrak{g}.$$ We consider $\langle \cdot, \cdot \rangle$ the Riemannian metric on $G$ given by left-translation of $\langle \cdot, \cdot \rangle_\mathfrak{g}$, \textit{i.e.} $\langle \cdot, \cdot \rangle_x = \langle (L_x)_\ast\cdot, (L_x)_\ast\cdot \rangle_\mathfrak{g}$, for every $x \in G$. Then the metric $\langle \cdot, \cdot \rangle$ is automatically left-invariant, and since $\langle \cdot, \cdot \rangle_{\mathfrak{g}}$ is Ad-invariant, we also have that $\langle \cdot, \cdot \rangle$ is right-invariant, meaning that both $(L_x)_\ast$ and $(R_x)_\ast$ are isometries of $\langle \cdot, \cdot \rangle$ for every $x \in G$. In order to define a Laplace-Beltrami operator on $G$, we need to fix a measure on $G$. So let $\mathrm{d}_G$ be the bi-invariant Haar meadure in $G$ (which exists because the group $G$ is compact), and consider $L^2(G)$ the $L^2$-space associated with this measure. Let $\Delta_G$ be the Laplace-Beltrami operator associated with this metric, \textit{i.e.} $\Delta_G = \mathrm{d}^\ast \mathrm{d}$, where $\mathrm{d}^\ast$ is the transpose of the de Rham operator with respect with the $L^2$-space, \textit{i.e.} $$\langle \mathrm{d} f, \omega \rangle_{L^2} = \langle f, \mathrm{d}^\ast \omega \rangle_{L^2},$$ for every $f \in \mathcal{C}^\infty(G)$ and $\omega \in \mathcal{C}^\infty(G, \mathrm{T}^\ast(G))$. Let $X_1, \dots, X_n$ be an orthonormal basis of $\mathfrak{g}$ with respect to $\langle \cdot, \cdot \rangle_{\mathfrak{g}}$, and let $\chi_1, \dots, \chi_n$ be the its dual basis. In view of $\mathrm{d}_G$ being bi-invariant, we have that $X_j^L$ and $X_j^R$ satisfies $$\langle X_j^L f, g \rangle_{L^2} = - \langle f, X_j^L g \rangle_{L^2},$$ and $$\langle X_j^R f, g \rangle_{L^2} = - \langle f, X_j^R g \rangle_{L^2},$$ for every $f,g \in \mathcal{C}^\infty(G)$. For each $j =1, \dots, n$ we set $\chi_j^L$ and the $\chi_j^R$ the following one forms:
    $$\chi_j^L(Y)(x) = \chi_j((L_{x^{-1}})_\ast Y),$$
    $$\chi_j^R(Y)(x) = \chi_j((R_{x^{-1}})_\ast Y),$$
    for every $x \in G$ and $Y \in \mathrm{T}_x G$, and here we are identifying $\mathfrak{g}$ with $\mathrm{T}_e G$. It follows from this definitions that $\chi_j^L$ is left-invariant and $\chi_j^R$ is right-invariant, \textit{i.e.} $(L_y)^\ast \chi_j^L = \chi_j^L$ and $(R_y)^\ast \chi_j$ for every $y \in G$. Since $\mathrm{d}_G$ and $\langle \cdot, \cdot \rangle$ are bi-invariant then we can write
    $$\mathrm{d}f = \sum_{j=1}^n (X_j^L f) \chi_j^L,$$
    and
    $$\mathrm{d}f = \sum_{j=1}^n (X_j^R f) \chi_j^R,$$
    for every $f \in \mathcal{C}^1(G)$. Therefore $\Delta_G = - \sum_{j=1}^n (X_j^L)^2 = -\sum_{j=1}^n (X_j^R)^2$. 
    \begin{lem}
    For every $X \in \mathfrak{X}_L(G)$ and $Y \in \mathfrak{X}_R(G)$ we have that $[\Delta_G, X] = [\Delta_G, Y] = 0$.
    \end{lem}
    \begin{proof}
    Let $Y \in \mathfrak{X}_R(G)$. Since for every $j=1, \dots, n$, $[Y, X_j^R] \in \mathfrak{X}_{R}(G)$, we can write $[Y, X_j^R] = \sum_{k=1}^n a_{j,k} X_j^R$, where $a_{j,k} = \langle [Y, X_j^R], X_k^R \rangle \in \mathbb{R}$, so
    \begin{align*}
        - [Y, \Delta_G] & = [Y, \sum_{j=1}^n (X_j^R)^2 ] \\
        & = \sum_{j=1}^n [Y, X_j^R]X_j^R + X_j^R[Y,X_j^R] \\
        & = \sum_{j,k = 1}^n a_{j,k} X_k^R X_j^R + a_{j,k} X_j^R X_k^R \\
        & = \sum_{j,k = 1}^n (a_{j,k} + a_{k,j}) X_k^R X_j^R.
    \end{align*}
    Since $\langle \cdot, \cdot \rangle$ is ad-invariant, we have that $$a_{j,k} = \langle [Y, X_j^R], X_k^R \rangle = -\langle [X_j^R, Y], X_k^R \rangle = - \langle X_j^R, [Y, X_k^R] \rangle = - a_{k,j},$$ so $[Y, \Delta_G] = 0$. The case for left-invariant vector fields is analogous.
    \end{proof}

    Let  $\hat{G}$ be the set of all classes of equivalences of irreducible, strongly continuous (s.c), unitary representations of $G$. Since $G$ is compact, the irreducible s.c. unitary representations of $G$ are finite dimensional. We recall that the Fourier transform of a distribution $u \in \mathcal{D}^\prime(G)$ is given by
    \[
        \mathfrak{F}_G[u](\xi) = \int_G u(x) \xi^\ast(x) \mathrm{d}_G(x), \qquad \xi \in \hat{G},
    \]
    where $\mathrm{d}_G$ stands for the bi-invariant Haar measure in $G$. We can recover the distribution $u$ with its Fourier transform in the following manner:
    \[
    u(x) = \sum_{\xi \in \hat{G}} d_\xi Tr\big(\mathfrak{F}_G[u](\xi) \xi\big),
    \]
    where $d_\xi$ is the dimension of $\xi$, and we recall that when $u$ is $L^2$ then the above series converges in the $L^2$ topology, and the same is valid replacing $L^2$ by $\mathcal{C^\infty}$, $\mathcal{C}^\omega$ and $G^s$. Furthermore, if $u \in L^2(G)$, then 
    \[
        \|u\|_{L^2}^2 = \sum_{\xi \in \hat{G}}d_\xi Tr\big(\mathfrak{F}_G[u](\xi)\mathfrak{F}_G[u](\xi)^\ast\big) = \sum_{\xi \in \hat{G}}d_\xi \|\mathfrak{F}_G[u](\xi)\|_{HS}^2. 
    \]
    The $\| \cdot \|_{HS}$ is the so-called Hilbert-Schmidt norm. More details about the $HS$-norm can be found in \cite{ruzhansky_pseudo-differential_2010}.
    
    Since $\langle \cdot, \cdot \rangle_\mathfrak{g}$ is Ad-invariant, for every $\xi \in \hat{G}$ there exists $\lambda_\xi$, such that
    \[
        \Delta \xi_{i,j} = \lambda_\xi \xi_{i,j}, \qquad i,j = 1, \dots, d_\xi.
    \]
    In view of $\Delta$ being a elliptic linear partial differential operator, by using standard arguments in the theory of linear PDEs, we have that, for every $k \in \mathbb{N}$, we can identify the Sobolev space $H^{2k}(G)$ with the space
    \[
        \left\{u \in \mathcal{D}^\prime(G) \; : \; \Delta^k u \in L^2(G) \right\},
    \]
    which is the same as
    \[
        \left\{u \in \mathcal{D}^\prime(G) \; : \; \sum_{\xi \in \hat{G}} d_\xi \langle \xi \rangle^{4k} \| \mathfrak{F}_G[u](\xi)\|_{HS}^2 < \infty \right\},
    \]
    where $\langle \xi \rangle = (1 + \lambda_\xi)^\frac{1}{2}$.  
    Since $H^{-2k}(G)$ is the topological dual of $H^{2k}(G)$, we can also identify this space with
    \[
        \left\{u \in \mathcal{D}^\prime(G) \; : \; \sum_{\xi \in \hat{G}} d_\xi \langle \xi \rangle^{-4k} \| \mathfrak{F}_G[u](\xi)\|_{HS}^2 < \infty \right\},
    \]
    \begin{dfn}
        We denote by $\mathcal{M}(\hat{G})$ the set of maps $F : \hat{G} \longrightarrow \cup_{m=1}^\infty U(\mathbb{C}^m)$ satisfying $F(\xi) \in U(\mathbb{C}^{d_\xi})$ for all $\xi$.
    \end{dfn}
    We note that the Fourier transform is an isometry between $H^{2k}(G)$, and the following weighted $L^2$ space:
    \[
        \left\{ F \in \mathcal{M}(\hat{G}) \; : \; \sum_{\xi \in \hat{G}} d_\xi \langle \xi \rangle^{4k} \|F(\xi)\|_{HS}^2 \right\}
    \]
    The same is valid for $H^{-2k}(G)$. We note that the Fourier transform in $\mathbb{R}^N$ is a isometry between $H^s(\mathbb{R}^N)$ and $L^2\big(\mathbb{R}^N; (1+|\xi|^2)^s\mathbb{d}\xi\big)$. Therefore, in view of Proposition \ref{prop:interpolation_L2} and Corollary \ref{cor:interpolation_isomorphism}, for every $s \in \mathbb{R}$, we have that $H^s(G)$ can be isometrically identified with
    \[
       \left\{u \in \mathcal{D}^\prime(G) \; : \; \sum_{\xi \in \hat{G}} d_\xi \langle \xi \rangle^{2s} \| \mathfrak{F}_G[u](\xi)\|_{HS}^2 < \infty \right\}, 
    \]
    with norm 
    \[
        \|u\|_s = \left( \sum_{\xi \in \hat{G}} d_\xi \langle \xi \rangle^{2s} \|\mathfrak{F}_G[u](\xi)\|_{HS}^2 \right)^\frac{1}{2}.
    \]
    
    \section{The local-global version}  \label{sec:local_global}

    Let $K$ be a compact Lie group. For $v \in \mathcal{E}^\prime (\mathbb{R}^n \times K)$, and $s_1, s_2 \in \mathbb{R}$, we define the following anisotropic Sobolev norms:
    \begin{equation}
        \| v \|_{s_1,s_2}^2 \doteq \sum_{\xi \in \hat{K}} \int_{\mathbb{R}^n} (1 + |\tau|^2)^{s_1} \langle \xi \rangle^{2s_2} \mathrm{d}_\xi \| \hat{v}(\tau, \xi) \|^2 \mathrm{d}\tau,
    \end{equation}
    where 
    \begin{align*}
        \hat{v}(\tau, \xi) & = \iint_{U \times K} v(t,x) e^{-i t \cdot \tau} \xi^\ast(x) \mathrm{d}\mu_K(x) \mathrm{d} t \\
        & = \int_K v^\tau(x)  \xi^\ast(x) \mathrm{d}\mu_K(x) \\
        & = \mathfrak{F}_K[v^\tau](\xi),
    \end{align*}
    and $\| \hat{v}(\tau, \xi) \|^2 = \|\mathfrak{F}_G[v^\tau](\xi) \|_{HS}$. With an analogous argument used in Section \ref{sec:preliminaries} using interpolation spaces and the the use of a finite covering of $K$ and a partition of unity, one can show that this anisotropic Sobolev norm is equivalent to the one coming from the local theory, thus we have that $\| v \|_{s_1,s_2} \leq \| v\|_{s_1 + s_2}$, where the last is the Sobolev norm on $\mathbb{R}^n \times K$, and we also note that $\| v \|_{0,0} = \| v\|_{0}$, the usual $L^2$ norm. 
    Let $U \subset \mathbb{R}^n$ be an open set. For every $s \in \mathbb{R}$, we say that $u \in H^\varepsilon_{\text{loc,glob}}(U \times K)$ if $\phi u \in H^s(\mathbb{R}^m \times K)$ for every $\phi \in \mathcal{C}^\infty_c(U)$.
    
    Let $P$ be a linear partial differential differential operator on $U \times K$ with smooth coefficients, let $t_0 \in U$ be fixed, and suppose that for every point $(t_0, x) \in U \times K$, there exists an $\varepsilon > 0$ such that $P$ is $\varepsilon$-subelliptic at $(0, x)$, \textit{i.e.} there exists $\tilde{U} \subset U$, an open neighborhood of $t_0$, and $V \subset K$, an open neighborhood of $x$, such that the following property holds true:
    \begin{equation*}
        u \in L^2_{\text{loc}}(\tilde{U} \times V), Pu \in L^2_{\text{loc}}(\tilde{U} \times V) \Rightarrow u \in H^\varepsilon_{\text{loc}}(\tilde{U} \times V).
    \end{equation*}
    Since $K$ is compact, there exist $\tilde{U} \subset U$, an open neighborhood of $t_0$, and $\varepsilon > 0$, such that 
    \begin{equation*}
        u \in L^2_{\text{loc,glob}}(\tilde{U} \times K), Pu \in L^2_{\text{loc,glob}}(\tilde{U} \times K) \Rightarrow u \in H^\varepsilon_{\text{loc,glob}}(\tilde{U} \times K).
    \end{equation*}
    By the Closed Graph Theorem, for every $\phi \in \mathcal{C}^\infty_c(\tilde{U})$ there exists a constant $C > 0$ such that
    \begin{equation}\label{eq:subelliptic_estimate_local_global}
        \| \phi u \|_\varepsilon^2 \leq C (\|\phi Pu\|_0^2 + \|\phi u\|_0^2), 
    \end{equation}
    for every $u \in L^2_{\text{loc,glob}}(\tilde{U} \times K)$. Using a standard technique, based on the Youngs' inequality for products, \textit{i.e.} $ab \leq \frac{1}{p}a^p + \frac{1}{q}b^q$, if $a,b>0$ and $\frac{1}{p} + \frac{1}{q} = 1$, one can prove that for every $k \in \mathbb{Z}_+$,
    \begin{align*}
        \|\phi u\|_0^2 & \leq \frac{\varepsilon}{\varepsilon + 2k}\left( \frac{4Ck}{\varepsilon + 2k}\right)^{\frac{2k}{\varepsilon + 2k}} \| \phi u\|_{0,-2k}^2 + \frac{1}{2C}\| \phi u\|_{0,\epsilon}^2 \\
        & \leq \frac{\varepsilon}{\varepsilon + 2k}\left( \frac{4Ck}{\varepsilon + 2k}\right)^{\frac{2k}{\varepsilon + 2k}} \| \phi u\|_{0,-2k}^2 + \frac{1}{2C}\| \phi u\|_\epsilon^2.
    \end{align*}
    So we can rewrite \eqref{eq:subelliptic_estimate_local_global} as
    \begin{equation*}
        \| \phi u \|_\varepsilon^2 \leq 2C \left(\|\phi Pu\|_0^2 + \frac{\varepsilon}{\varepsilon + 2k}\left( \frac{4Ck}{\varepsilon + 2k}\right)^{\frac{2k}{\varepsilon + 2k}} \| \phi u\|_{0,-2k}^2 \right).
    \end{equation*}
    Let $L_0, L_1, \ldots, L_\nu$ be as in \eqref{eq:vf2}, and $Q$ as in \eqref{eq:Q} . In view of condition (H) and Hörmander's Theorem, there exists $\tilde{U} \subset U$ an open neighborhood of the origin, and $C,\varepsilon > 0$ such that
    \begin{equation*}
        \| \phi u \|_\varepsilon^2 \leq 2C \left(\|\phi Qu\|_0^2 + \frac{\varepsilon}{\varepsilon + 2k}\left( \frac{4Ck}{\varepsilon + 2k}\right)^{\frac{2k}{\varepsilon + 2k}} \| \phi u\|_{0,-2k}^2 \right), 
    \end{equation*}
    for every $u \in\mathcal{C}^\infty(\tilde{U}\times K)$ and $k \in \mathbb{Z}_+$.
    
    Suppose that $Qu = f \in \mathcal{C}^\omega(U \times K)$. By Hörmander's Theorem, we have that $u$ is smooth. Let $\Delta_K$ be as in section \ref{sec:preliminaries}. Since $\phi$ is independent of $x$ (here we are writing $(t,x) \in U \times K$), and the metric considered on section \ref{sec:preliminaries} is ad-invariant, the Laplace-Beltrami operator $\Delta_K$ commutes with $Q - a$ (see section \ref{sec:preliminaries}) and with $\phi$. Therefore, we can apply the same technique as in the proof of Theorem 4.1 of \cite{braun_rodrigues_lower_2016} to obtain the following estimate
    \begin{equation}\label{eq:analytic_vector_Delta_K}
        \| \phi \Delta_K^k u \|_0 \leq B^{k+1}k!^2,
        \qquad \forall k \in \mathbb{Z}_+.
    \end{equation}
    Since the operator $\Delta_K$, viewed as an operator on $K$, is elliptic, we have that, on $U \times K$, its characteristic variety is given by
    \[
        \text{Char}(\Delta_K) = \{ (t,x,\tau,0) \in T^\ast(U\times K) \; : \; \tau \neq 0 \}. 
    \]
    It follows from \eqref{eq:analytic_vector_Delta_K} that $u$ is an analytic vector for $\Delta_K$, so by Theorem 3.2 of \cite{bolley_analyticite_1979} we have that $WF_a (u)|_{(t_0,x)} \subset \{(\tau, 0) \; : \; \tau \neq 0 \}$, for every $x \in K$. In view of condition (E2), we have that $Q$ is elliptic in the variable $t$, meaning that the points $\{ (t_0,x,\tau,0) \; : \; x \in K, \tau \in \mathbb{R}^n \}$ does not belong to $\text{Char}(Q)$, and since $WF_a (u) \subset \text{Char}(Q)$, we obtain that $WF_a (u)|_{(t_0,x)} = \emptyset$, implying that $\{t_0\} \times K \notin \text{singsupp}_a (u)$, and since $t_0 \in U$ is arbitrary, we have that $\text{singsupp}_a \,(u) = \emptyset$.

    \section{The reduction to the local-global version}\label{sec:reduction_local_global}

    By \cite[Corollary 10.1.11]{hilgert_structure_2012} and the fact that $\Omega$ is compact, there is a finite covering $\mathfrak U$ of $\Omega \doteq G / K$ by open sets such that for each $U \in \mathfrak U$, there is a smooth section $\sigma = \sigma_U$ for $\pi : G \to \Omega$ such that
    \begin{equation}
    \label{eq:hilgert_diffeo}
        (u, k) \in U \times K \mapsto \Psi(u,k) \doteq \sigma(u)k \in \sigma(U)K
    \end{equation} is a diffeomorphism onto an open subset $V = \sigma(U)K$ of $G$. Notice that we can take $\sigma$ to be real-analytic, and since the group multiplication is real-analytic, the map $\Psi$ is real-analytic. A similar argument also shows that $\Psi^{-1}$ is real-analytic.

    Let $g \in \Psi(U \times K)$. Writing $g = \Psi(u,k) = \sigma(u)k$,  we have that $u = \pi(g)$. Since $\sigma$ is a section of the projection map $\pi$ on $U$, \textit{i.e.} $\pi(\sigma(v)) = v$, for every $v \in U$, we have that, $\pi(g) = \pi(\sigma(u)) = u$, in other words, $g \sim \sigma(u)$, $g = \sigma(\pi(g))k$, for some $k \in K$, \textit{i.e.} $$k = \sigma(\pi(g))^{-1}g = \mathfrak{i}(\sigma(\pi(g)))g = \mu(\mathfrak{i}(\sigma(\pi(g))),g).$$ In conclusion, $\Psi^{-1}(g) = (\pi(g), \mu(\mathfrak{i}(\sigma(\pi(g)), g)) \in U \times K$. 
    %
    %
    %


    Let $\phi(g) \doteq i(\sigma \circ \pi(g)) g$ be defined for $g \in \pi^{-1}(U)$ and consider $g = \sigma(u)k$ and $g' = gk' = \sigma(u)kk'$ for some $k' \in K$.

    Notice that $\pi(g) = gK = gkk'K = \pi(g')$ and so we have that
    $$k = i(\sigma(\pi(g)))g = i(\sigma(\pi(g'))) g = i(\sigma(\pi(g'))) g'k'^{-1} = R_{k'^{-1}} \circ i(\sigma(\pi(g'))) g'$$
    which implies that $\phi(R_{k'}(g)) = R_{k'} \circ \phi(g)$ for all $g \in \pi^{-1}(U)$. This means that, if a vector field is invariant by the right-action of $K$ in $G$, then it descents into a right-invariant vector field in $K$ via $\phi_*$.
    


    Since $X_1, \dots, X_\nu$ satisfies \eqref{eq:e1}, the vector fields $L_j = \Psi^{-1}_*(X_j), j = 1, \ldots, \nu$ satisfy condition \eqref{eq:e2} and from basic properties of pushforward of vector fields, they also satisfy (H). 
    


    \section{Proof of Theorem \ref{thm:global}}

    Consider the vector fields $X_0, X_1, \ldots, X_\nu$ and the operator $P$ satisfying the conditions of Theorem \ref{thm:global}. We assume that $Pu \in C^\omega(G)$ and we want to show that $u \in C^\omega(G)$. In other words, we wish to prove that $\operatorname{singsupp}_a\,(u) = \emptyset$. In view of Hörmander's Theorem, we already know that $u \in \mathcal{C}^\infty(G)$. So let $p \in G$ be an arbitrary point. Then there exists $U$ a neighborhood of $\pi(p)$ such that $\Psi_U : U \times K \to \sigma(U)K = V$, as in \eqref{eq:hilgert_diffeo}, is a $\mathcal{C}^\omega$-diffeomorphism, and $p \in V$. As seen in section \ref{sec:reduction_local_global}, the pushfoward of the operator $P$ by $\Psi_U^{-1}$ is an operator $Q$ on $U\times K$ as \eqref{eq:Q} satisfying (H) and \eqref{eq:e2}. As seen on section \ref{sec:local_global}, this means that $u|_{V}\circ \Psi_U^{-1} \in \mathcal{C}^\omega(U\times K)$, and in particular, $\Psi_U^{-1}(p) \notin \text{singsupp}_a\,(u|_{V}\circ \Psi_U^{-1})$. Since the analytic singular support is invariant under $\mathcal{C}^\omega$-diffeomorphisms, we have that $p \notin \text{singsupp}_a\,(u)$. Since $p \in G$ is arbitrary, we have that $\text{singsupp}_a\,(u) = \emptyset$.
    
        


\section{Appendix}

    In this section, we recall some results and definitions regarding the theory of interpolation spaces and the K-Method. Let $\mathcal{E}$ be a topological vector space, and let $E_0, E_1$ be two normed spaces, continuously embedded in $\mathcal{E}$. We equip $E_0\cap E_1$ and $E_0 + E_1$ with the following norms
    \begin{align*}
        & \|a\|_{E_0\cap E_1} \doteq \max \{\|a\|_{E_0}, \|a\|_{E_1} \},\\
        & \|a\|_{E_0+E_1} \doteq \inf\{\|a_0\|_{E_0} + \|a_1\|_{E_1} \; : \; a = a_0 + a_1, a_0 \in E_0, a_1 \in E_1 \}.
    \end{align*}
    To simplify the notation, we write the infimum in the definition of $\|a\|_{E_0+E_1}$ as $\inf_{a = a_0 + a_1}\big(\|a_0\|_{E_0} + \|a_1\|_{E_1}\big)$.
    \begin{dfn}
        We say that a normed space $E$ is an intermediate space between $E_0$ and $E_1$ if $E_0 \cap E_1 \subset E \subset E_0 + E_1$, with continuous inclusions.
    \end{dfn}
    \begin{dfn}
        We say that a normed space $E$ is an interpolation space between $E_0$ and $E_1$ if it is an intermediate space and if it has the following property: if $T:E_0 + E_1 \longrightarrow E_0 + E_1$ is a linear map, and $T\in \mathcal{L}(E_0,E_0)\cap\mathcal{L}(E_1,E_1)$, then $T \in \mathcal{L}(E)$. The space $E$ is said to be of exponent $0 < \theta < 1$ if there exists a constant $C > 0$ such that 
        \[
            \|T\|_{\mathcal{L}(E,E)} \leq C\|T\|_{\mathcal{L}(E_0,E_0)}^{1-\theta}\|T\|_{\mathcal{L}(E_1,E_1)}^{\theta}, \qquad \forall T\in \mathcal{L}(E_0,E_0)\cap\mathcal{L}(E_1,E_1).
        \]
    \end{dfn}
    In the following, we describe a method to construct such interpolation spaces, it is called the K-Method. For every $a \in E_0 + E_1$ and every $t > 0$, we define
    \[
        K_2(t;a) \doteq \inf_{a = a_0 + a_1}\big(\|a_0\|_{E_0}^2 + t^2\|a_1\|_{E_1}^2 \big)^\frac{1}{2}.
    \]
    For every $0 < \theta < 1$, we define
    \[
        \big(E_0, E_1\big)_{\theta,2} \doteq \left\{ a \in E_0 + E_1 \; : \; t^{-\theta}K_2(t;a) \in L^2\left(\mathbb{R}_+;\frac{\mathrm{d}t}{t} \right) \right\},
    \]
    with norm 
    \[
        \|a\|_{(E_0,E_1)_{\theta,2}} \doteq \|t^{-\theta}K_2(t;a)\|_{L^2\left(\mathbb{R}_+;\frac{\mathrm{d}t}{t}\right)}.
    \]
    \begin{prp}
        Let $F_0$ and $F_1$ be normed spaces with the same properties as $E_0$ and $E_1$ before. If $T: E_0 + E_1 \longrightarrow F_0 + F_1$ is a linear map, such that $T \in \mathcal{L}(E_0,F_0) \cap \mathcal{L}(E_1,F_1)$, then $T \in \mathcal{L}((E_0,E_1)_{\theta,2}, (F_0,F_1)_{\theta,2})$ for every $0 < \theta < 1$, and $\|T\|_{\mathcal{L}((E_0,E_1)_{\theta,2}, (F_0,F_1)_{\theta,2})} \leq \|T\|_{\mathcal{L}(E_0,F_0)}^{1-\theta}\|T\|_{\mathcal{L}(E_1,F_1)}^{\theta}$.
    \end{prp}
    \begin{proof}
    Let $a \in E_0 + E_1$. For each decomposition $a = a_0 + a_1$, we have that
    \begin{align*}
        K_2(t;Ta) & \leq \big(\|T a_0\|_{F_0}^2 + t^2 \|T a_1\|_{F_1}^2\big)^{\frac{1}{2}} \\
        & \leq \big( \|T\|_{\mathcal{L}(E_0,F_0)}^2 \|a_0\|_{E_0}^2 + t^2\|T\|_{\mathcal{L}(E_1,F_1)}^2 \|a_1\|_{E_1}^2\big)^\frac{1}{2}\\
        & = \|T\|_{\mathcal{L}(E_0,F_0)} \Bigg(\|a_0\|_{E_0}^2 + t^2\frac{\|T\|_{\mathcal{L}(E_1,F_1)}^2}{\|T\|_{\mathcal{L}(E_0,F_0)}^2} \|a_1\|_{E_1}^2\Bigg)^\frac{1}{2}.
    \end{align*}
    Thus, $K_2(t;Ta) \leq \|T\|_{\mathcal{L}(E_0,F_0)} K_2\left(t\frac{\|T\|_{\mathcal{L}(E_1,F_1)}}{\|T\|_{\mathcal{L}(E_0,F_0)}}; a \right)$. To conclude, we note that
    \begin{align*}
        \|Ta\|_{(F_0,F_1)_{\theta,2}} & = \left(\int_0^\infty t^{-2\theta} K_2(t;Ta)^2 \frac{\mathrm{d}t}{t}\right)^{\frac{1}{2}} \\ 
        & \leq \left( \int_0^\infty t^{-2\theta}\|T\|_{\mathcal{L}(E_0,F_0)}^2 K_2\left(t\frac{\|T\|_{\mathcal{L}(E_1,F_1)}}{\|T\|_{\mathcal{L}(E_0,F_0)}}; a \right)^2 \frac{\mathrm{d}t}{t}\right)^\frac{1}{2} \\
        & \underset{t = \frac{\|T\|_{\mathcal{L}(E_0,F_0)}}{\|T\|_{\mathcal{L}(E_1,F_1)}}s}{=} \|T\|_{\mathcal{L}(E_0,F_0)}^{1-\theta}\|T\|_{\mathcal{L}(E_1,F_1)}^\theta \left( \int_0^\infty s^{-2\theta} K_2(s; a)^2 \frac{\mathrm{d}s}{s}\right)^\frac{1}{2} \\
        & = \|T\|_{\mathcal{L}(E_0,F_0)}^{1-\theta}\|T\|_{\mathcal{L}(E_1,F_1)}^\theta \|a\|_{(E_0,E_1)_{\theta,2}}. \qedhere
    \end{align*}
    \end{proof}
    \begin{cor}\label{cor:interpolation_isomorphism}
    Let $F_0$ and $F_1$ be normed spaces with the same properties as $E_0$ and $E_1$ before. If $T: E_0 + E_1 \longrightarrow F_0 + F_1$ is a linear map, such that $T$ is a linear isomorphism between $E_0$ and $F_0$, and between $E_1$ and $F_1$, then $T$ is a linear isomorphism between $\big(E_0,E_1\big)_{\theta,2}$ and $\big(F_0,F_1\big)_{\theta,2}$ for every $0 < \theta < 1$.
    \end{cor}
    We restrict our attention to the case in which $E_0$ and $E_1$ are weighted $L^2$ spaces. Let $\Omega \subset \mathbb{R}^N$ be an open set, and let $\omega$ be a positive measurable function defined on $\Omega$. We define
    \[
        E(\omega) \doteq \left\{ u: \Omega \longrightarrow \mathbb{C} \; : \; u\,\text{is measurable}\,, \int_\Omega |u(x)|^2 \omega(x) \mathrm{d}x < \infty \right\},
    \]
    and
    \[
        \|u\|_\omega \doteq \left( \int_\Omega |u(x)|^2 \omega(x) \mathrm{d}x \right)^\frac{1}{2}.
    \]

    With the introduced definitions and notation, we can easily prove the following proposition.
    
    \begin{prp}\label{prop:interpolation_L2}
        Let $\omega_0$ and $\omega_1$ be two positive measurable functions on $\Omega$. Then, for every $0 < \theta < 1$, 
        \[
        \big( E(\omega_0), E(\omega_1) \big)_{\theta,2} = E(\omega_\theta),
        \]
        where $\omega_\theta = \omega_0^{1-\theta} \omega_1^\theta$.
    \end{prp}
    \begin{proof}
    We have an explicit formula for $K_2(t;a) = \inf_{a = a_0 + a_1} \big( \|a_0\|_{\omega_0}^2 + t^2 \|a_1\|_{\omega_1}^2 \big)^\frac{1}{2}$. Let $a \in E(\omega_0) + E(\omega_1)$, and let $x \in \Omega$. Consider the following optimization problem: find $\lambda \in \mathbb{C}$ that minimizes 
    \[
        |\lambda|^2\omega_0(x) + t^2|a(x) - \lambda|^2\omega_1(x).
    \]
    Such lambda is obtained by solving the following equations:
    \begin{align*}
        &\lambda\omega_0(x) - t^2(a(x) - \lambda) \omega_1(x) = 0,\\
        &\overline{\lambda}\omega_0(x) - t^2\overline{(a(x) - \lambda)}\omega_1(x) = 0, 
    \end{align*}
    thus $\lambda = \frac{t^2 a(x)\omega_1(x)}{\omega_0(x) + t^2\omega_1(x)}$. Therefore
    \begin{align*}
        K_2(t;a)^2 &= \int_\Omega \frac{t^4\omega_1(x)^2 + t^2\omega_0(x)^2}{(\omega_0(x) + t^2\omega_1(x))^2}|a(x)|^2 \mathrm{d}x \\
        & = \int_\Omega \frac{t^2\omega_0(x)\omega_1(x)}{\omega_0(x) + t^2\omega_1(x)}|a(x)|^2 \mathrm{d}x,
    \end{align*}
    and this implies that
    \begin{align*}
        \|t^{-\theta}K_2(t;a)\|_{L^2\left(\mathbb{R}_+;\frac{\mathrm{d}t}{t}\right)}^2 &= \int_0^\infty \int_\Omega t^{-2\theta}\frac{t^2\omega_0(x)\omega_1(x)}{\omega_0(x) + t^2\omega_1(x)}|a(x)|^2 \mathrm{d}x \frac{\mathrm{d}t}{t} \\
        & = \int_\Omega \int_0^\infty t^{-2\theta}\frac{t^2\omega_0(x)\omega_1(x)}{\omega_0(x) + t^2\omega_1(x)}|a(x)|^2 \frac{\mathrm{d}t}{t} \mathrm{d}x  \\
        & \underset{t = \sqrt{\frac{\omega_0(x)}{\omega_1(x)}}s}{=} \int_\Omega \int_0^\infty \frac{s^{2-2\theta}}{1+s}\frac{\omega_0(x)\omega_1(x)}{\omega_0(x)^\theta\omega_1(x)^{1-\theta}}|a(x)|^2 \frac{\mathrm{d}s}{s} \mathrm{d}x\\
        & = \left(\int_0^\infty \frac{s^{1-2\theta}}{1+s} \mathrm{d}s\right) \|a\|_{\omega_\theta}^2. \qedhere
    \end{align*} 
    \end{proof}

\bibliographystyle{alpha}
\bibliography{references.bib, nicholas.bib}

\end{document}